\def\Z{{Z\!\!\!\! Z}}

\def\CR{\hbox{{$\cal R$}}}
\def\isom{{\cong}}
\def\tens{\mathop{\otimes}}
\def\la{{\triangleright}}
\def\id{{\rm id}}
\def\kfg{$k_FG$}

\documentstyle[11pt]{article}
\newtheorem{propos}{{\rm PROPOSITION}}
\newtheorem{corol}[propos]{{\rm COROLLARY}}
\newtheorem{defin}[propos]{{\rm DEFINITION}}

\begin{document}

{\ }\qquad  \vspace{1.2in}

\parindent 0pt

{\Large NEW APPROACH TO OCTONIONS AND CAYLEY ALGEBRAS}

\baselineskip 13pt
\bigskip
\bigskip
HELENA ALBUQUERQUE\footnote{Supported by CMUC-JNICT and by
Praxis 2/2.1/Mat7458/94}

\medskip
Departamento de Matematica-Faculdade de Ciencias e Tecnologia

Universidade de Coimbra, Apartado 3008

3000 Coimbra, Portugal

\bigskip
SHAHN MAJID\footnote{Royal Society University Research Fellow and
Fellow of Pembroke College, Cambridge}

\medskip
Department of Applied Mathematics and Theoretical Physics

University of Cambridge

Cambridge CB3 9EW, UK

\bigskip
\medskip ABSTRACT We announce a new approach to the
octonions as {\em quasiassociative} algebras. We strip out the
categorical and quasi-quantum group considerations in our longer
paper and present here (without proof) some of the more algebraic
conclusions.

\parindent 10pt

\section{{\Large INTRODUCTION}}

Usually one recognises the nonassociativity of the octonions by
saying that they are instead alternative algebras. While this is
true, the property of being alternative has a much weaker
character than associativity and, as a result, many standard ideas
and constructions for associative algebras do not go through in
the alternative case. In our paper \cite{AlbMa}, to which this
note is a short introduction, we have introduced a full solution
to this problem based on modern ideas from category theory and
Drinfeld's theory of quasiquantum groups.

Without going into any details (see \cite{AlbMa}), the new
formulation is that the octonions and other Cayley algebras live
naturally as objects in a monoidal category\cite{Mac}\cite{Ma}. For
any three objects $V,W,Z$ in such a category there is an associator
isomorphism $\Phi_{V,W,Z}:(V\tens W)\tens Z\to V\tens (W\tens Z)$
which performs the rebracketting. Mac Lane's pentagon condition on
$\Phi$ ensures that {\em we can do all constructions as if there
are no brackets} (i.e. as if $\tens$ is strictly associative).
After writing any desired constructions as the composition of
various maps, we simply insert $\Phi$ as needed for the
compositions to make sense, and all different ways to do this will
give the same result (this is Mac Lane's coherence theorem). So
working in such a category is no harder than usual associative
linear algebra. For example, an algebra $A$ in such a category
means
\[\bullet\circ(\bullet\tens\id)=\bullet\circ(\id\tens\bullet)
\circ\Phi_{A,A,A}\]
for the product $\bullet$, where $\Phi$ is inserted for the
bracketting to make sense. So, recognizing the octonions as such a
{\em quasiassociative algebra} (or {\em quasialgebra} for short)
makes them as good as associative in the precise sense explained
above.

In \cite{AlbMa} we introduce and study a class of such
quasialgebras that contain composition algebras and more general
algebras obtained by a generalised Cayley-Dickson process. All
algebras are considered over a field $k$ of characteristic
different from 2. The required class of quasialgebras arises
naturally by a certain `Drinfeld twisting'\cite{Dri}\cite{Ma} or
deformation of classical group algebras, as follows.

\section{{\Large QUASIALGEBRAS $k_FG$}}

First of all, we know that if we consider the set of complex
numbers, the quaternions or the octonions, all these algebras have
something in common: If we choose a suitable basis and remove the
$\pm$ signs from the multiplication tables of these algebras, we
have the tables of the additive groups $G=\Z_2$ (for complex
numbers), $G=\Z_2\times \Z_2$ (for quaternions) and $G=\Z_2\times
\Z_2 \times \Z_2$ (for octonions)\cite{Dix}. We view the signs in the
multiplication tables as an invertible 2-cochain $F:G\times G\to
k$ (a nowhere-zero function which is 1 when either argument is the
group identity $e\in G$). Writing $F(x,y)=(-1)^{f(x,y)}$, one has
explicitly\cite{AlbMa},
\[  G=\Z_2,\quad f(x,y)=xy,\quad x,y\in\Z_2\quad\quad\quad\quad
{\rm (Complex\ numbers)},\]
\[  G=(\Z_2)^2,\quad
f(\vec{x},\vec{y})=x_1y_1+(x_1+x_2)y_2\quad\quad\quad
{\rm (Quaternions)},\]
where $\vec{x}=(x_1,x_2)\in G$ is a vector notation
and the components $x_1,x_2$ are viewed in the field $\Z_2$.
\[ G=(\Z_2)^3,\quad
f(\vec{x},\vec{y})= \sum_{i\le
j}x_iy_j+y_1x_2x_3+x_1y_2x_3+x_1x_2y_3\quad {\rm (Octonions)},\]
where $\vec{x}=(x_1,x_2,x_3)\in G$ is a vector
notation. Similarly for higher Cayley algebras where $G=(\Z_2)^n$.

{}From the group $G$ and the cochain $F$ we recover the complex,
quaternion and octonion algebras as the `deformation' $k_FG$ of the
group algebra of the appropriate $G$. This is the vector space with
basis labeled by $G$ and the product
\[ x\bullet y=F(x,y) xy\]
for $x,y\in G$, where $xy$ is the group product in $G$. So this is
a kind of deformation of the usual group algebra of $G$. In
quantum groups we do the deformation by introducing a parameter
$q$ such that when $q$ tends to 1 we have the original algebras.
Here we do the deformation by introducing a cochain $F$.

\begin{propos}\cite{AlbMa} Let $G$ be a group and $F$ any
invertible 2-cochain. Then $k_FG$ is a $G$-graded quasialgebra with
associator $\Phi$ determined by the coboundary $\phi$ of $F$.
Explicitly, it is quasiassociative in the sense \[ (x\bullet
y)\bullet z=\phi(x,y,z)x\bullet (y\bullet z)\] for all $x,y,z\in
G$, and
\[ \phi(x,y,z)= {{F(x,y)F(xy,z)}\over {F(y,z)F(x,yz)}}.\]
\end{propos}

To explain the setting here, for $G$ a group and $\phi:G\times
G\times G\to k$ an invertible group 3-cocycle, the category of
$G$-graded vector spaces becomes monoidal with the associator
$\Phi$ determined by the 3-cocycle $\phi$ and the grading. A
quasialgebra with this form of $\Phi$ is called a {\em $G$-graded
quasialgebra}. It consists of an algebra $A$, a $G$-grading which
respects the product and unit (so the degree of $1\in A$ is $e\in
G$, the group identity), and the quasiassociativity law
 \[ (a\bullet b)\bullet c=\phi(|a|,|b|,|c|)a\bullet (b\bullet c)\]
for all elements $a,b,c\in A$ of degree $|a|,|b|,|c|$. In our
case, $k_FG$ is such a $G$-graded quasialgebra with $\phi$ built
from $F$ and with $|x|=x$ for $x\in G$.

For the complex number and the quaternion algebras, $F$ is closed,
i.e. $\phi$ is trivial and the algebras happen to be strictly
associative. This is because $f$ in these cases is quadratic. As
soon as we introduce a cubic or higher `interaction' term in $f$,
as in the case of the octonions, $\phi$ typically becomes
nontrivial and the algebra $k_FG$ nonassociative. In the case of the
Octonions it is
\[ \phi(\vec{x},\vec{y},\vec{z})=(-1)^{(\vec{x}\times\vec{y})
\cdot\vec{z}}\]
(the vector cross product and vector dot product in the exponent, 
i.e. the determinant $|\vec{x}\ \vec{y}\ \vec{z}|$). On the other
hand, if we simply drop the cubic or higher terms in the above
family, we clearly obtain the corresponding Clifford algebra with
negative signature\footnote{We would like to thank Tony Smith for
asking us to clarify this point} (the relations are immediate and
the dimensions match) i.e. these are obtained as $k_FG$ with
\[ G=(\Z_2)^n,\quad f(\vec{x},\vec{y})=\sum_{i\le j}x_iy_j
\quad\quad\quad\quad {\rm (Clifford\ algebras)},\]
as the associative version of the octonion or Cayley algebra, which
is another way to see the close relationship between these and
Clifford algebras. The positive signature algebras are obtained
similarly with $i<j$ in $f$.

Also, we are mainly interested in $G$ Abelian and specialise to
this case from now on. For $\phi$ of the coboundary form, the
category of $G$-graded spaces is symmetric\cite{Mac}, i.e. for any
two objects $V,W$ there is a generalised transposition isomorphism
$\Psi_{V,W}:V\tens W\to W\tens V$. A quasialgebra $A$ is {\em
quasicommutative} if $\bullet=\bullet\circ\Psi_{A,A}$. This is the
case for all $k_FG$ with $\Psi$ determined by a function
$\CR(x,y)=F(x,y)/F(y,x)$. Explicitly,
\[ x\bullet y=\CR(x,y) y\bullet x,\]
for all $x,y\in G$. For complex numbers, quaternions and octonions,
etc., the function $\CR$ has the simple form
\[ \CR(x,y)=\cases{1&if\ $x=e$\ or\ $y=e$\ or\ $x=y$\cr
-1&otherwise.}\]
We call this important case {\em altercommutative}. By contrast,
for the Clifford algebras one has
\[ \CR(\vec{x},\vec{y})=(-1)^{\sum_{i\ne j} x_iy_j}\]
which is not in general altercommutative (for $n>2$).

\section{{\Large ALGEBRAIC PROPERTIES OF $k_F G$}}

We have just seen that the functions $\phi,\CR$ built from $F$
allow for the categorical setting of the algebras $k_FG$ as
quasiassociative and quasicommutative. In particular, the algebra
is associative iff $\phi=1$ and commutative iff $\CR=1$. We now
summarise less obvious results expressing the more conventional
algebraic properties of $k_FG$ in terms of these functions
$\phi,\CR$. We refer to \cite{AlbMa} for proofs and details.

\begin{propos}\cite{AlbMa} \kfg\ is an alternative algebra
{\em iff}
\[ \phi^{-1}(y,x,z) + \CR(x,y)\phi^{-1}(x,y,z)=1+\CR(x,y)\]
\[ \phi(x ,y  ,z   )+{\CR(z,y)}\phi(x ,z   ,y  )=1+{\CR(z,y)}\]
for all $x,y,z\in G$. In this case,
\[  \phi(x,x,y)=\phi(x,y,y)=\phi(x,y,x)=1 \] for all
$x,y\in G$.
\end{propos}
\vskip2mm

Next we consider involutions. Since we have a special basis of
$k_FG$ it is natural to consider involutions diagonal in this
basis.
\vskip2mm
\begin{propos}\cite{AlbMa} \kfg\ admits an involution which is
diagonal in the basis $G$ {\em iff} $\CR(x,y)={{s(x)s(y)}\over
s(xy)}$ for some 1-cochain $s:G\to k$ (a nowhere-zero function with
$s(e)=1$) obeying $s^2=1$. In this case, one has
$\CR(x,y)=\CR(y,x)$ and $\phi(x,y,z)=\phi(z,y,x)^{-1}$ for all
$x,y,z\in G$.
\end{propos}

The corresponding involution here is $\sigma(x)=s(x)x$ for all
$x\in G$. Let $A$ be a finite dimensional algebra with identity
element 1 and let $\sigma$ be an involution in $A$. We say that
$\sigma$ is a strong involution if $a+\sigma(a)$,
$a\bullet\sigma(a)\in k1$ for all $a\in A$.

\begin{propos}\cite{AlbMa} $k_FG$ admits a diagonal strong
involution $\sigma$ {\em iff}

i) $G\simeq (\Z_2)^n$ for some $n$,

ii) $\sigma(e)=e$, $\sigma(x)=-x$ for all $x\ne e$,

iii) $k_FG$ is altercommutative.
\end{propos}

Given $k_FG$ we have a natural function $s(x)=F(x,x)$ and consider
now the possibility of defining a strong involution using this. For
all statements of simplicity in the following we assume $|G|>2$.

\begin{propos}\cite{AlbMa} If $\sigma(x)=F(x,x)x$ for all
$x\in G$ defines a strong involution, then the algebra \kfg\ is
simple and the following are equivalent

i) $k_FG$ is an alternative algebra,

ii) $k_FG$ is a composition algebra.
\end{propos}

Finally, it is possible to characterize a natural class of $k_F G$
algebras that are composition algebras,

\begin{propos}\cite{AlbMa} Let $k_FG$ admit a strong diagonal
involution $\sigma (x)= s(x)x$. Then $q(x)=x\bullet
\sigma (x)$ makes $k_FG$ a composition algebra {\em iff}

i) $s(xy)F(x,y)^2 F(xy,xy)=s(x)s(y)F(x,x)F(y,y)$, for all $x,y\in
G$.

ii) $F(x,xz)F(y,yz)F(z,z)s(z)+F(x,yz)F(y,xz)F(xyz,xyz)s(xyz)=0$,
for all $x,y\in G$ with $x\not=y$.
\end{propos}

An important corollary of the last result is:

\begin{corol}\cite{AlbMa} If $G\isom(\Z_2)^n$ then the Euclidean
norm quadratic function defined by $q(x)=1$ for all $x\in G$ makes
$k_FG$ a composition algebra {\em iff}

i) $F^2(x,y)=1$ for all $x,y\in G$

ii) $F(x,xz)F(y,yz)+F(x,yz)F(y,xz)=0$, for all $x,y,z\in G$ with
$x\ne y$.
\end{corol}

In this case $\sigma(x)=F(x,x)x$ for all $x\in G$ is a strong
involution and $k_FG$ is simple and alternative.

\section{{\Large CAYLEY-DICKSON PROCESS FOR $k_F G$}}

We have a generalisation of the Cayley-Dickson process as follows.
Again, details are in \cite{AlbMa}.

\begin{defin}
Let $G$ be an Abelian group $F$ a 2-cochain on it.
 For any 1-cochain $s:G\to k$ and $\alpha\ne 0$ we define
 $\bar G=G\times \Z_2$ and on it the 2-cochain $\bar F$ and
1-cochain $\bar s$,
\[\bar F(x,y)=F(x,y),\quad \bar
F(x,vy)=s(x)F(x,y),\quad \bar F(vx,y)=F(y,x),\]
\[\bar F(vx,vy)=\alpha s(x)F(y,x),\quad \bar s(x)=s(x),\quad \bar
s(vx)=-1\]
 for all $x,y\in G$. Here $x\equiv(x,e)$ and
$vx\equiv(x,\nu)$ denote elements of $\bar G$, where
$\Z_2=\{e,\nu\}$ with product $\nu^2=e$.
\end{defin}

We say that $k_{\bar F}\bar G$ is the {\em generalised
Cayley-Dickson extension} of $k_FG$ associated to $s,\alpha$. The
motivation is that if $\sigma(x)=s(x)x$ is a strong involution,
then $k_{\bar F}\bar G$ is the usual Cayley-Dickson extension of
$k_FG$ associated to $\sigma,\alpha$. Note that since all unital
composition algebras over $k$ are obtained by repeated
Cayley-Dickson extension\cite{ZSSS}, they are all of the form of a
quasialgebra $k_F G$ in the last proposition of the preceding
section with $G$ a power of $\Z_2$.

The natural application for our generalised Cayley-Dickson process
is when $k_FG$ admits a diagonal involution $\sigma(x)=s(x) x$
(but not necessarily strong). In this case we have:

\begin{propos}\cite{AlbMa} The 3-cocycle $\bar\phi$ of
$k_{\bar F}\bar G$ is given by
\[\bar \phi     (x ,y  ,z   )=\phi(x ,y,z ),\quad
\bar \phi (v x ,y ,z )={\CR(y,z)}\phi(x ,y ,z ),\]
\[ \bar\phi (x ,v y ,z )={\CR(y,z)}{\CR(xy,z)}
\phi(x ,y ,z),\quad
\bar\phi(x ,y ,v z )={\CR(x,y)}\phi(x ,y ,z ),\]
\[\bar \phi      (v x  ,v y ,z )={\CR(xy,z)}\phi(x   ,y,z ),
\quad\bar \phi   (v x   ,y  ,v z )={\CR(y,z)}{\CR(x,y)}
\phi(x ,y ,z),\]
\[ \bar \phi (x ,v y ,v z )=\CR(x,yz)\phi(x ,y ,z ),\quad
\bar \phi (v x     ,v y       ,v z )={\CR(xy,z)}{\CR(x,y)}\phi(
x,y ,z ),\] for $x,y,z\in G$.
\end{propos}

Using this calculation, one may show under the same assumptions:

\begin{propos}\cite{AlbMa}

i) $k_{\bar F}\bar G$ is associative {\em iff} $k_FG$ is
associative and commutative.

ii) If $k_FG$ has trivial centre then $k_{\bar F}\bar G$ is
alternative {\em iff} $k_FG$ is associative and $s(x)=-1$ for all
$x\in G$ and $x\ne e$.

iii) If $s$ defines a strong diagonal involution then $k_{\bar F}\bar
G$ is alternative {\em iff} $k_FG$ is associative.
\end{propos}

This extends to more general $k_FG$ some well-known considerations
for octonions and higher Cayley algebras.

\section{{\Large CONCLUDING REMARKS}}

We conclude with some remarks about other classes if
quasialgebras. In fact, the input data for our $k_FG$ construction
is clearly very general. If we denote the elements of the finite
group $G$ by $\{x_1=e,x_2,\dots,x_n\}$ then we can represent the
cochain $F$ by an $n\times n$ matrix with entries
$F_{ij}=F(x_i,x_j)$. It has $1$ in the first row and column and
all entries non-zero. Conversely, any such matrix will do for a
cochain and yield a quasialgebra. Therefore it is a wide-open
question what other groups and cochains might be interesting; here
we list just a few natural classes.

First of all, motivated by composition algebras for the Euclidean
norm (isomorphic to complex, quaternions or octonions), where the
cochain is represented by certain normalized Hadamard matrices, a
natural more general class of examples is $k_FG$ where $F$ is a
general normalised Hadamard matrix. Some results for (and
low-dimensional examples of) this kind of $k_FG$ are given in
\cite{AlbMa}. Hadamard matrices have even dimension but odd
dimensional examples of $k_FG$ can be obtained by taking $F$ with
first column and row $1$ and the rest an Hadamard matrix.

Another general choice of $F$, overlapping with the extended
Hadamard case, is (for $G$ any group),
\[ F(x,y)=\cases{1&if\ $x=e$\ or\ $y=e$\ or $x\ne y$\cr
-1&otherwise},\]
where one may show that $k_FG$ is simple, commutative and in
general nonassociative.

Other interesting examples of $k_FG$ quasialgebras come from the
theory of finite fields, and include the generalisations of
octonions based on Galois sequences in \cite{Dix}.

Finally, we would also like to recall that our approach answers
such questions as what means a representation of the octonions.
Following the method explained in the introduction, a
representation of the quasialgebra $k_FG$ means a $G$-graded
vector space $V$ and a degree-preserving action $\la$ obeying
\[ (x\bullet y)\la v=\phi(x,y,|v|)
x\la(y\la v)\] for all $v\in V$. This is explained in
\cite{AlbMa}, where it is also shown that a representation is
equivalent to an algebra map from $k_FG$ to a certain quasialgebra
${\rm End}_\phi(V)$ of {\em quasimatrices} associated to any $V$.

\end{document}